\documentclass[12pt]{article}
\usepackage[utf8]{inputenc}
\usepackage{amssymb,amsthm,bm,graphicx,xcolor}
\usepackage{mathtools}
\usepackage[T1]{fontenc}

\newtheorem{thm}{Theorem}[section]

\newtheorem{cor}[thm]{Corollary}
\newtheorem*{cor*}{Corollary}
\newtheorem{lem}[thm]{Lemma}

\theoremstyle{definition}

\newtheorem{defi}[thm]{Definition}

\newtheorem{cla}[thm]{Claim}
\newtheorem{rek}[thm]{Remark}

\newcommand{\N}{\mathbb{N}}

\newcommand{\R}{\mathcal{R}}

\newcommand{\bp}{\begin{proof}}
\newcommand{\ep}{\end{proof}}
\newcommand{\eq}[1]{\begin{align*}#1\end{align*}}

\newcommand{\rb}{\{R,B\}}
\newcommand{\dcube}{2^{[2n]}}
\newcommand{\cube}{2^{[n]}}
\renewcommand{\ss}{S\subseteq [2n]}
\newcommand{\ceilh}{\ceil{n/2}}
\newcommand{\floorh}{\floor{n/2}}

\DeclarePairedDelimiter\ceil{\lceil}{\rceil}
\DeclarePairedDelimiter\floor{\lfloor}{\rfloor}

\title{A construction for Boolean cube Ramsey numbers}
\author{Tom Bohman\thanks{Department of Mathemathics, Carnegie Mellon University. Correspondence to tbohman@math.cmu.edu. This work was supported by a grant from the Simons
Foundation (587088, TB)} \, Fei Peng\thanks{Department of Mathematics, National University of Singapore.}}
\date{\today}

\begin{document}

\maketitle
\begin{abstract}
Let $Q_n$ be the poset that consists of all subsets of a fixed $n$-element set, ordered by set inclusion.
The poset cube Ramsey number $R(Q_n,Q_n)$ is defined as the least~$m$ such that any 2-coloring of the elements of $Q_m$ admits a monochromatic copy of $Q_n$. The trivial lower bound $R(Q_n,Q_n)\ge 2n$ was improved by Cox and Stolee, who showed $R(Q_n,Q_n)\ge 2n+1$ for $3\le n\le 8$ and $n\ge 13$ using a probabilistic existence proof. In this paper, we provide an explicit construction that establishes $R(Q_n,Q_n)\ge 2n+1$ for all $n\ge 3$. The best known upper bound, due to Lu and Thompson, is $ R(Q_n, Q_n) \le n^2 - 2n + 2$.
\end{abstract}
\section{Introduction}
A central theme of combinatorics is the fact that large discrete systems often contain subsystems with a higher degree of organization than the original system. Results of this flavor appear in many areas. One example is the Erd\H{o}s-Szekeres Theorem, which states that every sequence of $ab+1$ real numbers has a monotonously increasing subsequence of length $a$ or a monotonously decreasing subsequence of length $b$.
Another example is the Erd\H{o}s-Szekeres Conjecture, which asks for the minimum number $n$ such that any $n$ points in general position in the plane contain $k$ points in convex position \cite{suk}. Graph Ramsey theory \cite{graham1990ramsey} studies the graph-theoretic analog of this theme. Ramsey's Theorem for graphs states that for all $k$, there exists $n$ such that any edge 2-coloring of a clique (complete graph) of size $n$ contains a monochromatic clique of size $k$. The least such $n$ is denoted $R(k,k)$.

In this paper, we will focus on the Ramsey problem for poset cubes. A partially ordered set, or a \textit{poset}, is a set equipped with a partial relation $\le$ that is reflexive, anti-symmetric and transitive. The cube $Q_n$, which is the power set of $[n]:=\{1,2, \dots, n\}$ equipped with the inclusion relation, plays an important role in the theory of posets, and we consider the natural Ramsey question in this context. We use the power set $2^{[n]}$ itself to refer to the associated poset when it is unambiguous. A \textit{poset embedding} is an order-preserving (and nonorder-preserving) injection $f$ from a poset $\mathcal{P}$ to another poset $\mathcal{P}'$: that is, $x\le_\mathcal{P}y\iff f(x)\le_{\mathcal{P}'}f(y)$. In the context of cubes, a function $f:2^{[n]}\to2^{[m]}$ is an embedding if it is injective and satisfies $S\subseteq T\iff f(S)\subseteq f(T)$. In this case, we say that the range of $f$ -- which is a subset of $2^{[m]}$ -- is a \textit{copy} of $2^{[n]}$.

This paper focuses on the following Ramsey theoretic quantity introduced by Axenovich and Walzer \cite{axenovich2017boolean}. Given $n\in\N$, define the (poset) \textit{cube Ramsey number} $R(Q_n,Q_n)$ to be the least integer $m$ such that every 2-coloring of the elements of $Q_m$ admits a monochromatic copy of $Q_n$. Note that the individual subsets of $[m]$ are colored, instead of the inclusion relations between them; in other words, we consider a vertex-coloring instead of an edge-coloring. 

It is not hard to see that $R(Q_n,Q_n)\ge 2n$: that is, one can 2-color the cube $Q_{2n-1}$ without any monochromatic copy of $Q_n$. The method is to simply color all odd-sized subsets red and even-sized subsets blue. In fact, if we group subsets by their cardinality, we would split $Q_{2n-1}$ into $2n$ \textit{layers}, and any division of the layers into $n$ totally-red layers and $n$ totally-blue layers will suffice.

The best lower bound on $R(Q_n,Q_n)$ to date is slightly better than $2n$. Cox and Stolee \cite{cox2018ramsey} established that $R(Q_n,Q_n)\ge 2n+1$ for $3\le n\le 8$ and $n\ge 13$. Their argument for $n\ge 13$ relies on a family of special sets, which is shown to exist using the Lovász Local Lemma, and there is no known explicit construction of this family. Here, we exhibit an explicit 2-coloring of $Q_{2n}$ without a monochromatic copy of  $Q_n$ for all $n\ge 4$ which, when combined with the previously known construction for $n=3$ \cite{axenovich2017boolean}, gives:
\begin{thm}\label{theorem:main}
$\forall n \ge 3, R(Q_n,Q_n)\ge 2n+1.$
\end{thm}
\noindent
A generalization of Theorem~\ref{theorem:main} for large $n$ was recently established by Gr\'osz, Methuku and Tompkins\cite{new}. The best known upper bound of the cube Ramsey number is $R(Q_n,Q_n)\le n^2-2n+2$ established in \cite{lu2019poset}, and a recent paper \cite{falgas2020existence} contains some results regarding cube Ramsey numbers of a random poset, and also obtains $R(Q_3,Q_3)=7$. 

We mention in passing the closely-related poset variant of the classical Turán problem, which asks how many elements we can choose from $2^{[n]}$ without seeing a copy of $2^{[k]}$. This is unsolved even when $k=2$ \cite{axenovich2012q} (that is, picking elements to not see a ``diamond''), and many believe that the answer is ${n\choose \floorh}+{n\choose \floor{(n-1)/2}}$, witnessed by taking the two layers closest to the middle. The best known upper bound is $(( \sqrt{2} + 3)/2  + o(1)) {n\choose \floor{n/2}}$  \cite{diamondupper}.

\section{The coloring}
We exhibit a 2-coloring of $\dcube$ admitting no monochromatic copy of $\cube$. To define the coloring, we first group $[2n]$ into $n$ pairs:
\begin{defi}
$a,b\in [2n]$ are said to be in the same \textit{pair} if $\ceil {a/2} = \ceil{b/2}$. So the pairs are $\{1,2\},\{3,4\}\dots$. We say a set $S$ \textit{has} a pair if $S$ contains that pair, and $S$ \textit{misses} the pair if $S$ contains neither element of the pair. The \textit{partner} of $a \in [2n]$ is the unique $b$ such that $ \{a,b\}$ is a pair. An element $x \in S$ is a \textit{single} in $S$ if the partner of $x$ is not an element of $S$.
\end{defi}
We sometimes characterize a set by its number of pairs and singles: for example, $\{1,2,3,5,6,7,8\}$ is said to have 3 pairs and 1 single. 

\begin{defi}
A collection $\R \subseteq \dcube$ is \textit{pair-enforcing} if for every $\ss$ s.t. $\ceilh\le |S|<n$, $S \in\R$ only if $S$ has a pair. 
\end{defi}

\begin{defi}
A collection $\R \subseteq \dcube$ is \textit{miss-forbiding} if for every $\ss$ s.t. $n<|S|\le n+\floorh$, $S \in\R$ only if $S$ doesn't miss any pair.
\end{defi}

\begin{defi}
A collection $\R \subseteq \dcube$ is \textit{not-too-high} if for every $S \in\R$ we have $|S|\le n+\floorh$.
\end{defi}

\begin{defi}
A collection $\R \subseteq \dcube$ is \textit{flip-susceptible} if for all $S_1,S_2\subseteq [2n]$ such that $|S_1|=|S_2|=n$, $|S_1\cup S_2|=n+1$ and neither $S_1$ nor $S_2$ has any pair, at most one of $S_1$ and $S_2$ is in $\R$.
\end{defi}
It is named as such because $S_1$ and $S_2$ would be all the same except for a choice-flip about one specific pair.

\begin{defi}
A collection $\R \subseteq \dcube $ is \textit{restrictive} if it is pair-enforcing, miss-forbiding, not-too-high and flip-susceptible.
\end{defi}

Now we are ready to present the construction.

\begin{defi}\label{definition:coloring}
Define the coloring $c_0:\dcube\to \rb$ as follows: $\forall \ss,$
\begin{itemize}
    \item If $|S|<\ceilh$, $S$ is red.
    \item If $\ceilh\le |S| < n$, $S$ is red iff $S$ has a pair.
    \item If $|S|=n$, $S$ is red iff the sum of its elements is odd.
    \item If $n<|S|\le n+\floorh$, $S$ is red iff $S$ does not miss any pair.
    \item  If $|S|>n+\floorh$, $S$ is blue.
\end{itemize}
\end{defi}

Let $\R$ be the collection of subsets of $ [2n]$ that are colored red by $ c_0$. Note that $\R$ is restrictive. 
Note further that the collection of sets of the form $ \bar{S}:= [2n] \setminus S$ such that $S$ is colored 
blue by $ c_0$ is also restrictive. Indeed, if $\ceilh\le |S| < n$ and $S \not\in \R$ then $S$ contains 
no pair and therefore $ \bar{S} $ does not miss any pair, and if $n<|S|\le n+\floorh$ and $ S \not\in \R$ then $S$ 
misses a pair and therefore $ \bar{S}$ contains a pair.
Note further that if $ \mathcal{Q} \subset \dcube$ is a copy of $\cube$ then the collection $ \{ \bar{S}: S \in \mathcal{Q} \}$ is also a copy of $\cube$. So, as \cite{axenovich2017boolean} gives a construction that establishes $ R( Q_3,Q_3) \ge 7$, 
Theorem \ref{theorem:main} follows from the following result:

\begin{thm}\label{theorem:restrictivenocopy}
If $n\ge 4$, a restrictive collection $ \R \subseteq \dcube$ does not contain a copy of $\cube$.
\end{thm}
We prove this in the next Section.


%
%

\section{Proof of Theorem \ref{theorem:restrictivenocopy}}

We begin with two preliminary observations about poset embeddings. For the next few results and notations, we fix a poset embedding $f:2^X\to2^Y$, where $X$ and $Y$ are finite sets.

\begin{lem}\label{lemma:postcard}
If $A \subseteq B\subseteq X$ then $|f(B)|-|f(A)|\ge |B|-|A|$. In particular (setting $A=\emptyset$), $|f(B)|\ge |B|$.
\end{lem}
\begin{proof}
Say $|B|-|A|=k$. Then there exists a chain $A=S_0\subsetneq S_1\subsetneq \dots \subsetneq S_k=B$. Since $f$ is an embedding, $f(A)=f(S_0)\subsetneq f(S_1)\subsetneq \dots\subsetneq f(S_k)=f(B)$. Hence, $|f(B)|-|f(A)|\ge k$.
\end{proof}

\begin{defi}
With respect to the poset embedding $f:2^X\to 2^Y$, the \textit{top element} refers to $f(X)$, and the \textit{top children} refers to all immediate children of the top element, namely all sets of the form $f(X\setminus \{a\})$ ($a\in X$).
\end{defi}
\begin{lem}\label{lemma:topch}
If the top element has cardinality $N$, then the intersection of any $k$ top children ($1\le k\le |X|$) has cardinality at most $N-k$; that is, $\forall\text{ nonempty } I\subseteq X$, \eq{
\bigg|\bigcap_{a\in I} f(X\setminus \{a\})\bigg|\ \le\ N-|I|.
}
\end{lem}
\bp
First note that every top child has size at most $N-1$, so the lemma is true when $k=1$. Then, $\forall \text{ nonempty } I\subseteq X, \forall j\in X\setminus I$, \eq{\bigcap_{a\in I} f(X\setminus \{a\})\ \supsetneq\ \bigcap_{a\in I\cup \{j\}} f(X\setminus \{a\})} because $f(\{j\})$ is a subset of the former but not the latter. Thus, when $I$ grows from a singleton to a $k$-element set, the size of the intersection decreases to at most $N-k$.
\ep

With these observations in hand, we are ready to consider restrictive collections.

\begin{lem}
Let $\R$ be a restrictive collection. For all $S \in \R$, there is some $S^+\supseteq S,\ |S^+|=n+\floorh$ such that $\R \cup \{S^+\} $ is a restrictive collection.
\end{lem}
\bp
Let $ S \in \R$ and
assume for the sake of contradiction that
$S$ has more than $\floorh$ pairs. Then $|S| > n$. Since $\R$ is miss-forbidding and not-too-high, $S\in \R$ only if $S$ doesn't miss any pair. Note that a set that misses no pair and has $k$ pairs has cardinality exactly $n+k$, so we have $|S|>n+\floorh$, witnessing $S\notin \R$. Hence, the assumption is false.

As $S$ has no more than $\floorh$ pairs, by adding elements as needed, one can find $S^+\supseteq S$ that has exactly $\floorh$ pairs and doesn't miss any pair. Such $S^+$ will have size $n+\floorh$. The collection $R\cup\{S^+\}$ contains $S^+$ while continuing to be restrictive.
\ep

\begin{defi}
A copy of $\cube$ is \textit{maximal} if its top element has size $n+\floorh$.
\end{defi}

\begin{cor}\label{corollary:tomax}
If a restrictive collection $\R$ contains a copy of $\cube$, then some restrictive collection $\R^+\supseteq \R$ contains a maximal copy of $\cube$.
\end{cor}

By Corollary \ref{corollary:tomax}, it suffices to show that there does not exist a restrictive collection that contains a maximal copy of $ \cube.$

In the remainder of this Section, we fix a restrictive collection $\R \subseteq \dcube$ and assume for the sake of contradiction that $\R$ contains a maximal copy of $\cube$. Let this maximal copy be the image of the embedding $ f : 2^{[n]} \to 2^{[2n]}$. 

By definition, the top element of this maximal copy has size $n+\floorh$. By the miss-forbidding property, this element cannot miss any pair, so it has exactly $\floorh$ pairs and $n-\floorh = \ceilh$ singles. Denote the pairs as $p_1,\dots, p_{\floorh}$ and the singles $s_1,\dots, s_{\ceilh}$. 
For ease of notation we let $\Pi$ denote the collection of pairs $p_1,\dots,p_{\floorh}$, and we let $\Sigma$ denote the collection of singles 
$s_1,\dots,s_{\ceilh}$.

\begin{cla}\label{lemma:allpairs}
In the maximal copy, if a top child misses an element of $\Sigma$ 
then the top child has size at most $n$ and has all $\floorh$ pairs in $\Pi$.
\end{cla}

\bp
Consider such a top child $S$ that misses an element from $\Sigma$.  This top child misses the pair that single lies in, because even the top element doesn't have the partner. Since $\R$ is miss-forbidding and not-too-high, $|S|\le n$.

As $S$ is a top child, we can write $S=f([n]\setminus\{a\})$ for some $a\in [n]$.
Sets in the copy that are below or equal to $S$ naturally form a copy of $2^{[n]\setminus \{a\}}$, with $S$ being the top element. 
There are $n-1$ top children in the smaller cube; by Lemma \ref{lemma:topch}, their intersection has size at most $n-(n-1)=1$. 

Suppose $S$ has at most $\floorh - 1$ of the pairs in $\Pi$. For each pair it has, one of its immediate children doesn't have it -- this is because at most one element of $ [2n]$ is in the intersection of all immediate children of $S$. Thus, there is a collection of at most $\floorh - 1$ immediate children of $S$ such that for every pair in $\Pi$ that $S$ has, one immediate child in the collection does not have that pair. 

These immediate children are of the form $f([n]\setminus\{a,b_i\})$, where $i$ ranges from 1 to $k$, and $k\le \floorh - 1$. Consider the set $S_1=f([n]\setminus \{a,b_1,\dots,b_k\})$. By Lemma \ref{lemma:postcard}, its size is at least $|[n]\setminus \{a,b_1,\dots,b_k\}|=n-(k+1)\ge \ceilh$, and since it is below $S$, its size is smaller than $n$. But on the other hand, $S_1$ is a subset of every $f([n]\setminus \{a,b_i\})$, so it does not have any pair. As $\R$ is pair-enforcing, $S_1$ cannot be in $\R$. But $S_1$ is in the copy. To resolve the contradiction, $S$ must have all the $\floorh$ pairs.
\ep

\begin{cla}\label{corollary:allsingle}
The maximal copy has at most $\ceilh-\floorh+1$ top children (i.e., one if $n$ is even and two if $n$ is odd) which miss at least one element of $\Sigma$.
\end{cla}

\bp
If $n$ is even, this is clear: by Claim \ref{lemma:allpairs}, such a top child has size at most $n$ but has all $\floorh = n/2$ pairs in $\Pi$, so it must consist of exactly the pairs in $\Pi$ and nothing else. Thus there is at most one such child.

If $n$ is odd, the pairs are determined for the same reason, but there could be an extra single in the top child. Assume for the sake of contradiction that three such top children exist. Since they are mutually incomparable, they all have an extra single from $\Sigma$. Say the top children are \eq{f([n]\setminus\{a\})&=p_1\sqcup\dots\sqcup p_{\floorh}\sqcup \{ s_i\},\\f([n]\setminus\{b\})&=p_1\sqcup\dots\sqcup p_{\floorh}\sqcup \{s_j\},\\f([n]\setminus\{c\})&=p_1\sqcup\dots\sqcup p_{\floorh}\sqcup \{s_k\},}
where $ \sqcup$ denotes disjoint union.

Then, \eq{f(\{a\})&\subseteq f(([n]\setminus\{b\})\cap ([n]\setminus\{c\}))\\&\subseteq f([n]\setminus\{b\})\cap  f([n]\setminus\{c\})\\&=p_1\sqcup\dots\sqcup p_{\floorh}\\&\subseteq f([n]\setminus\{a\}).}
This is a contradiction.
\ep


\begin{cla}\label{lemma:seq}
If $n\ge 4$ then there are $b_1,\dots,b_{\floorh}$ and $\hat{b}$, distinct elements of $[n]$, such that $f([n]\setminus\{b_i\})$ doesn't have $p_i$ for all $i\in[\floorh]$, and $f([n]\setminus\{\hat{b}\})$ doesn't have $p_j$ for some $j\in[\floorh]$.
\end{cla}

\bp
By Claim \ref{corollary:allsingle}, at least $n-(\ceilh-\floorh + 1)$ top children have all the elements of $\Sigma$. Construct a bipartite graph with bipartition $\Pi\sqcup \mathcal{C}$ where $\mathcal{C}$ is the collection of the top children that have all singles in $\Sigma$, and $p\in \Pi,\ S\in \mathcal{C}$ are connected iff $S$ does \textit{not} have the pair $p$. 

We claim that $\forall P\subseteq \Pi, |N(P)|\ge |P|$, where $ N(P)$ is the neighborhood of $p$ in this bipartite graph. If $P$ consists of $x>0$ pairs and $|N(P)|\le x-1$, then at least $n-(\ceilh-\floorh + 1)-(x-1)$ top children have all the $x$ pairs in $P$. By Lemma \ref{lemma:topch}, their intersection has size at most $(n+\floorh)-(n-(\ceilh-\floorh + 1)-(x-1))=\ceilh + x$. But $s_1,\dots, s_{\ceilh}$ and the $x$ pairs in $P$ are already in the intersection, which implies that the intersection has size at least $\ceilh+2x$. Because $x>0$, this is impossible.

By Hall's Theorem, there is a $\Pi$-saturating matching. That means, there are $b_1,\dots,b_{\floorh}$ all distinct such that for every $i$, $f([n]\setminus \{b_i\})$ doesn't have $p_i$. 

Finally, since $n-(\ceilh-\floorh + 1)>\floorh$ when $n\ge 4$, let $f([n]\setminus\{\hat{b}\})$ be an unused top child in $\mathcal{C}$. Since it has all the singles, it cannot have all the pairs (or it would be the top element), so there is some $p_j$ it does not have.
\ep
We are now ready to complete the proof of Theorem \ref{theorem:restrictivenocopy}.
Consider $b_1,\dots, b_{\floorh}$ and $\hat{b}$ from Claim \ref{lemma:seq}. Say $\hat{b}$ doesn't have $p_j$ and consider the sets \eq{S_1 
&\ =\ f([n]\setminus \{b_1,\dots,b_{j-1},b_j,b_{j+1},\dots,b_{\floorh}\}),\\S_2&\ =\ f([n]\setminus \{b_1,\dots,b_{j-1},\hat{b},~b_{j+1},\dots,b_{\floorh}\}),\\S_\lor &\ =\ f([n]\setminus \{b_1,\dots,b_{j-1},~~~~b_{j+1},\dots,b_{\floorh}\}),\\S_\cup &\ =\ S_1\cup S_2.}

Note that $S_1,S_2,S_\lor \in \R$. Because $S_1$ lies in every $f([n]\setminus \{b_i\})$, it doesn't have any pair. Since $\R$ is pair-enforcing, to have $S_1\in\R$ we must have either $|S_1|<\ceilh$ or $|S_1|\ge n$. Because $|S_1|\ge |[n]\setminus \{b_1,\dots,b_{\floorh}\}|=\ceilh$, we have $|S_1|\ge n$. Analogously, $|S_2|\ge n$. 
Since $S_1\neq S_2$, $|S_\cup|\ge n+1$.

On the other hand, by Lemma \ref{lemma:postcard},  \eq{|f([n])| - |S_\lor|\ \ge\  |\{b_1,\dots,b_{j-1},b_{j+1},\dots,b_{\floorh}\}|\ =\ \floorh-1.} Thus, \eq{|S_\lor|&\ \le\ |f([n])| -(\floorh - 1)\\&\ =\ (n+\floorh)-(\floorh -1)\\&\ =\ n+1.}

Note that $S_1,S_2\subseteq S_\lor$, so $S_\cup\subseteq S_\lor$. Thus, $|S_\cup|=n+1$ and $|S_1|=|S_2|=n$. Note again that neither $S_1$ nor $S_2$ has any pair. Because $\R$ is flip-susceptible, at most one of $S_1$ and $S_2$ is in $\R$. This is a contradiction. Hence, no restrictive collection can contain a maximal copy of $\cube$. In view of Corollary \ref{corollary:tomax}, the proof is complete.

\begin{rek}
The condition $ n \ge 4$ in Theorem~\ref{theorem:restrictivenocopy} is necessary. Indeed,
the construction we present here does have monochromatic copies of $ 2^{[3]}= 2^{[n]}$ in the case $n=3$. It turns out that the construction for the $n=3$ case given in \cite{axenovich2017boolean} is similar to our construction as the only difference is in the coloring of sets $S$ with $ |S|=n=3$.

Note that the Red-Blue colorings of the elements of $ Q_{2n}$ with the property that the
Red color class is flip susceptible while the complement of the Blue color class is also
flip susceptible are fixed on
sets of size $n$ which do not have a pair (and therefore have one element from
every pair). Indeed, we can view the collection of 
such sets as a graph by joining two such sets with an edge if their union has cardinality $n+1$, this graph
is bipartite, and in order to have the flip susceptible conditions
the partite sets in this bipartite graph must be monochromatic. Note that while our result requires fixed 
colors for all sets that do not have cardinality $n$ it places no condition on 
the colors of sets of size $n$ that contain a pair; these sets could be colored arbitrarily. 
The construction of \cite{axenovich2017boolean} 
succeeds in the case $n=3$ by placing conditions on these colors. It is tempting 
to think that this flexibility in the coloring is an indication that the lower bound on $ R( Q_n, Q_n)$ 
that we present here can be improved.

\end{rek}


\vskip2mm

\noindent {\bf Acknowledgment.} We thank the anonymous referees for helpful comments.

\vskip2mm

\noindent {\bf Data Availability Statement:} All data 
generated or analyzed during this study are included in this 
published article 

\bibliographystyle{alpha}
\bibliography{bib}

\begin{thebibliography}{FRMTZ20}

\bibitem[AMM12]{axenovich2012q}
Maria Axenovich, Jacob Manske, and Ryan Martin.
\newblock {$Q_2$}-free families in the {Boolean} lattice.
\newblock {\em Order}, 29(1):177--191, 2012.

\bibitem[AW17]{axenovich2017boolean}
Maria Axenovich and Stefan Walzer.
\newblock Boolean lattices: Ramsey properties and embeddings.
\newblock {\em Order}, 34(2):287--298, 2017.

\bibitem[CS18]{cox2018ramsey}
Christopher Cox and Derrick Stolee.
\newblock Ramsey numbers for partially-ordered sets.
\newblock {\em Order}, 35(3):557--579, 2018.

\bibitem[FRMTZ20]{falgas2020existence}
Victor Falgas-Ravry, Klas Markstr{\"o}m, Andrew Treglown, and Yi~Zhao.
\newblock Existence thresholds and {Ramsey} properties of random posets.
\newblock {\em Random Structures \& Algorithms}, 57(4):1097--1133, 2020.

\bibitem[GMT18]{diamondupper}
D\'aniel Gr\'osz, Abhishek Methuku, and Casey Tompkins.
\newblock An upper bound on the size of diamond-free families of sets.
\newblock {\em J. Combin. Theory Ser. A}, 156:164--194, 2018.

\bibitem[GMT21]{new}
D\'aniel Gr\'osz, Abhishek Methuku, and Casey Thompkins.
\newblock Ramsey numbers of boolean lattices.
\newblock {\em arXiv preprint. arXiv:2104.02002}, 2021.

\bibitem[GRS90]{graham1990ramsey}
Ronald~L Graham, Bruce~L Rothschild, and Joel~H Spencer.
\newblock {\em Ramsey theory}, volume~20.
\newblock John Wiley \& Sons, 1990.

\bibitem[LT22]{lu2019poset}
Linyuan Lu and Joshua~C Thompson.
\newblock Poset {Ramsey} numbers for {Boolean} lattices.
\newblock {\em Order}, 39:171--185, 2022.

\bibitem[Suk17]{suk}
Andrew Suk.
\newblock On the {Erdős}-{Szekeres} convex polygon problem.
\newblock {\em J. Amer. Math. Soc.}, 30(4):1047--1053, 2017.

\end{thebibliography}

\vskip2mm

\noindent
{\bf Correspondence to: tbohman@math.cmu.edu}

\end{document}